\documentclass{amsart}

\usepackage{amscd,amssymb,amsmath,graphicx,verbatim}
\usepackage[dvips]{hyperref}
\usepackage{textcomp}
\usepackage[OT1,T1]{fontenc}

\newtheorem{theorem}{Theorem}[section]

\theoremstyle{definition}
\newtheorem{definition}[theorem]{Definition}

\theoremstyle{remark}
\newtheorem{remark}[theorem]{Remark}

\numberwithin{equation}{section}

\begin{document}

\title[A HIDDEN SYMMETRY RELATED TO THE RH]{A HIDDEN SYMMETRY RELATED TO THE RIEMANN HYPOTHESIS WITH
THE PRIMES INTO THE CRITICAL STRIP}
\author{Stefano Beltraminelli}
\address{S. Beltraminelli, CERFIM, Research Center for Mathematics
and Physics, PO Box 1132, 6600 Locarno, Switzerland}
\email{stefano.beltraminelli@ti.ch}
\author{Danilo Merlini}
\address{D. Merlini, CERFIM, Research Center for Mathematics and
Physics, PO Box 1132, 6600 Locarno, Switzerland}
\email{merlini@cerfim.ch}
\author{Sergey Sekatskii}
\address{S. Sekatskii, Laboratoire de Physique de la Mati\`ere Vivante,
IPMC, BSP, Ecole Polytechnique F\'ed\'erale de Lausanne, 1015 Lausanne,
Switzerland}
\address{S. Sekatskii, Institute of Spectroscopy, Russian Academy
of Sciences, Troitsk Moscow region, Russia}
\email{serguei.sekatski@epfl.ch}
\label{I1}\label{I2}\label{I3}
\date{14.3.2008}
\begin{abstract}
In this note concerning integrals involving the logarithm of the
Riemann Zeta function, we extend\ \ some treatments given in previous
pioneering works on the subject and introduce a more general set
of Lorentz measures. We first obtain two new equivalent formulations
of the Riemann Hypothesis (RH). Then with a special choice of the
measure we formulate the RH as a ``hidden symmetry'', a global symmetry
which connects the region outside the critical strip with that inside
the critical strip. The Zeta function with all the primes appears
as argument of the Zeta function in the critical strip. We then
illustrate the treatment by a simple numerical experiment. The representation
we obtain go a little more in the direction to believe that RH may
eventually be true.
\end{abstract}
\maketitle

\section{Some integrals involving the Zeta function}

We start with some integrals involving the absolute values of the
logarithm of the Zeta function. As far as we know, the first work
in this direction is due to Wang, who discovered a RH criterium
involving these integrals \cite{28}. More recent pioneering works
are due to Volchkov \cite{22} who found an integral relation on
the complex plane with two variables equivalent to the Riemann Hypothesis
(RH). Later Balazard, Saias and Yor \cite{23}, established another
equivalence to the RH by an integral involving only one variable
i.e. by integration on the critical line. In a subsequent treatment
by one of us the analytical computations were extended to every
line perpendicular to the {\itshape x} axis to obtain an equivalence
to the RH involving explicitly $\mathfrak{R}( s) =\rho $, with the
appearance of a shift along the real axis of exactly $\frac{1}{2}$
\cite{24}.

In the present note we first extend some of the above mentioned
treatments by introducing a more general Lorentz measure which we
are free to normalize in order to obtain more simple formulas i.e.:
\[
C\operatorname*{\int }\limits_{-\infty }^{\infty }\frac{1}{{\rho
_{0}}^{2}+t^{2}}dt=C\frac{\pi }{\rho _{0}}=1
\]

\noindent so that $C=\frac{\rho _{0}}{\pi }$.

Then following the same calculation as in \cite{24}, we may establish
the following starting formula (equivalent to the RH, too) which
reads:
\begin{equation}
\varphi =\frac{\rho _{0}}{\pi }\operatorname*{\int }\limits_{-\infty
}^{\infty }\frac{\ln ( \left| \zeta ( \rho +i t) \right| ) }{{\rho
_{0}}^{2}+t^{2}}dt=\ln \left( \frac{\zeta ( \rho +\rho _{0}) \left( \rho
+\rho _{0}-1\right) }{\left| 1-\rho \right| +\rho _{0}}\right) %
\label{XRef-Equation-29183221}
\end{equation}

\noindent where $\rho >\frac{1}{2}$ and $\rho _{0}$ is supposed
to be a positive function of $\rho$ in particular an absolute value.
We remark, (\ref{XRef-Equation-29183221}) has been calculated supposing
RH is true and moreover it is equivalent to RH. (\ref{XRef-Equation-29183221})
contains a shift of amount $\rho _{0}$ and is more rich then the
previous established formulas since many choices of the function
$\rho _{0}$ are possible and a special choice of it may be more
convenient for the calculations.

Then since (see for example \cite{25}):
\[
\frac{\rho _{0}}{\pi }\operatorname*{\int }\limits_{-\infty }^{\infty
}\frac{\ln ( \left| 1-\rho +i t\right| ) }{{\rho _{0}}^{2}+t^{2}}dt=\ln
( \left| 1-\rho \right| +\rho _{0}) 
\]

\noindent we have
\begin{equation}
\varphi _{1}=\frac{\rho _{0}}{\pi }\operatorname*{\int }\limits_{-\infty
}^{\infty }\frac{\ln ( \left| \zeta ( \rho +i t) \left( \left| 1-\rho
\right| +i t\right) \right| ) }{{\rho _{0}}^{2}+t^{2}}dt=\ln ( \zeta
( \rho +\rho _{0}) \left( \rho +\rho _{0}-1\right) ) %
\label{XRef-Equation-29192731}
\end{equation}

For the special choice $\rho +\rho _{0}=1$, in the critical strip
($\frac{1}{2}<\rho <1$) we then obtain $\varphi _{1}=0$. Indeed,
this statement holds also for $\rho =\frac{1}{2}$ when it is nothing
else than the Balazard-Saias-Yor equality $\operatorname*{\int }\limits_{-\infty
}^{\infty }\frac{\ln ( |\zeta ( \frac{1}{2}+i t) |) }{\frac{1}{4}+t^{2}}dt=0$
\cite{23}. On the other hand:
\[
\frac{d}{d \rho _{0}}\varphi _{1}=\frac{1}{\rho _{0}}\varphi _{1}-2\rho
_{0}\frac{\rho _{0}}{\pi }\operatorname*{\int }\limits_{-\infty
}^{\infty }\frac{\ln ( \left| \zeta ( \rho +i t) \left( \left| 1-\rho
\right| +i t\right) \right| ) }{{\left( {\rho _{0}}^{2}+t^{2}\right)
}^{2}}dt=\frac{\zeta ^{\prime }( \rho +\rho _{0}) }{\zeta ( \rho
+\rho _{0}) }+\frac{1}{\rho +\rho _{0}-1}
\]

\noindent and we define now:
\begin{equation}
\varphi _{2}=\frac{\rho _{0}}{\pi }\operatorname*{\int }\limits_{-\infty
}^{\infty }\frac{\ln ( \left| \zeta ( \rho +i t) \left( \left| 1-\rho
\right| +i t\right) \right| ) }{{\left( {\rho _{0}}^{2}+t^{2}\right)
}^{2}}dt%
\label{XRef-Equation-29191719}
\end{equation}

We now make a more specific choice and set $\rho =\frac{1}{2}+\alpha
\text{}$ and $\rho +\rho _{0}=1$ in the critical strip, thus $\rho
_{0}=\frac{1}{2}-\alpha $ ($0<\alpha <\frac{1}{2}$ in the critical
strip). Outside the critical strip we may choose $\rho _{0}=\alpha
-\frac{1}{2}$, with $\alpha >\frac{1}{2}$. $\varphi _{1}$ and $\varphi
_{2}$ written below as a function of $\alpha$ constitute the first
Theorem of this note. 
\begin{theorem}

\label{XRef-Theorem-2121326}
\begin{equation}
\begin{aligned}
\varphi _{1}&=\frac{\left| -\frac{1}{2}+\alpha \right| }{\pi }\operatorname*{\int
}\limits_{-\infty }^{\infty }\frac{\ln ( \left| \zeta ( \frac{1}{2}+\alpha
+i t) \left( \left| -\frac{1}{2}+\alpha \right| +i t\right) \right|
) }{{\left( \frac{1}{2}-\alpha \right) }^{2}+t^{2}}dt \\ &=\left\{ \begin{array}{lcl}
 0 & , & 0<\alpha <\frac{1}{2} \\
 \ln ( \zeta ( 2\alpha ) \left( 2\alpha -1\right) )  & , & \alpha
>\frac{1}{2}
\end{array}\right. %
\label{XRef-Equation-213225656}
\end{aligned}
\end{equation}
\begin{equation}
\begin{aligned}
\varphi _{2}&=\frac{\left| -\frac{1}{2}+\alpha \right| }{\pi }\operatorname*{\int
}\limits_{-\infty }^{\infty }\frac{\ln ( \left| \zeta ( \frac{1}{2}+\alpha
+i t) \left( \left| -\frac{1}{2}+\alpha \right| +i t\right) \right|
) }{{\left( {\left( \frac{1}{2}-\alpha \right) }^{2}+t^{2}\right)
}^{2}}dt \\ &=\left\{ \begin{array}{lcl}
 -\frac{\gamma }{2\left( \frac{1}{2}-\alpha \right) } & , & 0<\alpha
<\frac{1}{2} \\
 -\frac{1}{2\alpha -1}\left( \frac{\zeta ^{\prime }( 2\alpha ) }{\zeta
( 2\alpha ) }+\frac{1}{2\alpha -1}\right) +\frac{\ln ( \zeta ( 2\alpha
) \left( 2\alpha -1\right) ) }{2{\left( \alpha -\frac{1}{2}\right)
}^{2}} & , & \alpha >\frac{1}{2}
\end{array}\right. %
\label{XRef-Equation-21322575}
\end{aligned}
\end{equation}

where $\gamma$ is the Euler-Mascheroni constant. 
\end{theorem}

The above formulas are equivalent to the RH. With the choice we
have considered, $\varphi _{1}$ and $\varphi _{2}$ both diverge
at $\alpha =\frac{1}{2}$ i.e. at the right border of the critical
strip. We note that $\varphi _{1}$ may be seen as a potential.

Numerical computations concerning (\ref{XRef-Equation-213225656})
and (\ref{XRef-Equation-21322575}) of Theorem \ref{XRef-Theorem-2121326}
may be done and presented as an illustration but from the known
numerical results on the zeros, the two functions are exact up to
${10}^{-20}$ in the critical strip and exact outside the critical
strip. So we omit here a numerical computation which will be presented
below for another case concerning the ``hidden symmetry'' which
we now introduce.
\begin{remark}

From (\ref{XRef-Equation-29183221}) it is easily seen that another
simple choices of $\rho _{0}$ give rise to a potential $\varphi$
of (\ref{XRef-Equation-29183221}) without any divergence. We simply
mention the case $\rho _{0}( \rho ) =\rho $ where (\ref{XRef-Equation-29183221})
gives:
\[
\varphi =\frac{\rho }{\pi }\operatorname*{\int }\limits_{-\infty
}^{\infty }\frac{\ln ( \left| \zeta ( \rho +i t) \right| ) }{\rho
^{2}+t^{2}}dt=\left\{ \begin{array}{lll}
 \ln ( \zeta ( 2\rho ) \left( 2\rho -1\right) )  & , & \frac{1}{2}<\rho
<1 \\
 \ln  \left( \zeta ( 2\rho ) \right)  & , & \rho >1
\end{array}\right. 
\]
\end{remark}

We then continue to establish another theorem which on the RH expresses
a kind of ``hidden symmetry''. To do this it is more convenient
to consider the potential $\varphi$ above instead of $\varphi _{1}$.
$\varphi$ is in fact a function of $\alpha$ which is not injective
and we may ask for what ($\alpha$, $\alpha ^{\prime }$) the potential
is the same that is $\varphi ( \alpha ) =\varphi ( \alpha ^{\prime
}) $.

\section{A ``hidden symmetry"}

We reconsider (\ref{XRef-Equation-29183221}) and set (in the critical
strip as before) $\rho =\frac{1}{2}+\alpha $ and $\rho _{0}=\frac{1}{2}-\alpha
$, $0<\alpha <\frac{1}{2}$. Then from (\ref{XRef-Equation-29183221})
we have for $ 0<\alpha <\frac{1}{2}$:
\begin{equation}
\varphi ( \alpha ) =\ln \left( \frac{1}{1-2\alpha }\right) %
\label{XRef-Equation-210235740}
\end{equation}

A very simple formula for the potential inside the critical strip.
For the potential outside (\ref{XRef-Equation-29183221}) gives:
\begin{equation}
\varphi ( \alpha ^{\prime }) =\ln ( \zeta ( \rho +\rho _{0}) ) =\ln
( \zeta ( 2\alpha ^{\prime }) ) %
\label{XRef-Equation-210235848}
\end{equation}

\noindent where we have set $\rho +\rho _{0}=2\alpha ^{\prime }$,
$\alpha ^{\prime }>\frac{1}{2}$.

We now observe that $\varphi ( \alpha ) $ of (\ref{XRef-Equation-210235740})
is increasing with $\alpha$ while $\varphi ( \alpha ^{\prime })
$ from (\ref{XRef-Equation-210235848}) is decreasing with $\alpha
^{\prime }$. Thus $\varphi$ is not injective in the interval [$\frac{1}{2},\infty
$[; this suggest the following definition.
\begin{definition}

The ``hidden symmetry'' is defined by the solution $(\alpha ,\alpha
^{\prime })$ of the equation:
\begin{equation}
\varphi ( \alpha ) =\ln \left( \frac{1}{1-2\alpha }\right) =\varphi ( \alpha
^{\prime }) =\ln ( \zeta ( 2\alpha ^{\prime }) ) 
\end{equation}
\end{definition}

\noindent where we write for the unique solution $\alpha ^{\prime
}=g( \alpha ) $.

We are still free to define for any $0<\alpha <\frac{1}{2}$ the
map:
\[
\alpha \longrightarrow 2g( \alpha ) +\alpha -\frac{1}{2}=\rho (
\alpha ) 
\]

\noindent and thus $\rho _{0}( \alpha ) =\frac{1}{2}-\alpha $ outside
the critical strip, so that $2\alpha ^{\prime }=\rho +\rho _{0}=2g(
\alpha ) $. Moreover:
\begin{equation}
\alpha =\frac{1}{2}\left( 1-\frac{1}{\zeta ( 2\alpha ^{\prime })
}\right) =\frac{1}{2}\left( 1-\prod \limits_{p\,\mathrm{prime}}\left(
1-\frac{1}{p^{2\alpha ^{\prime }}}\right) \right) %
\label{XRef-Equation-21201724}
\end{equation}

We may also write:
\begin{equation}
\alpha =\frac{1}{2}\left( 1-\frac{1}{\zeta ( 2\alpha ^{\prime })
}\right) =\frac{1}{2}\left( 1-\sum \limits_{n=1}^{\infty }\frac{\mu
( n) }{n^{2\alpha ^{\prime }}}\right) 
\end{equation}

\noindent where $\mu ( n) $ is the the M\"obius function of argument
{\itshape n}.\ \ 

Returning to (\ref{XRef-Equation-29183221}) we have that:
\[
\varphi ( \alpha ) =\frac{2\rho _{0}}{\pi }\operatorname*{\int }\limits_{0}^{\infty
}\frac{\ln ( \left| \zeta ( \rho +i t) \right| ) }{{\rho _{0}}^{2}+t^{2}}dt
\]

\noindent with the variable change $z=1-\frac{2\rho _{0}}{\rho _{0}+i
t}$ which maps the line $\rho _{0}+i t$ in the unit circle i.e.
$z= e ^{-2i \arctan ( \frac{t}{\rho _{0}}) }= e ^{-2i \theta }$,
where $\frac{d}{dt}\theta =\frac{1}{\rho _{0}( 1+\frac{t^{2}}{{\rho
_{0}}^{2}}) }$, we obtain:
\begin{equation}
\begin{aligned}
\varphi ( \alpha ) &=\frac{2}{\pi }\operatorname*{\int }\limits_{0}^{\pi
/2}\ln ( \left| \zeta ( \rho +i \rho _{0}\tan  \theta ) \right|
) d\theta \\ &=\frac{2}{\pi }\operatorname*{\int }\limits_{0}^{\pi /2}\ln
\left( \left| \zeta ( \frac{1}{2}+\alpha +i \left( \frac{1}{2}-\alpha
\right) \tan  \theta ) \right| \right) d\theta %
\label{XRef-Equation-2120580}
\end{aligned}
\end{equation}

\noindent  inside the critical strip. Outside the critical strip,
using the choice written above we obtain:
\begin{equation}
\begin{aligned}
\varphi ( \alpha ^{\prime }) &=\frac{2}{\pi }\operatorname*{\int
}\limits_{0}^{\pi /2}\ln ( \left| \zeta ( \rho +i \rho _{0}\tan
\theta ) \right| ) d\theta \\ &=\frac{2}{\pi }\operatorname*{\int }\limits_{0}^{\pi
/2}\ln \left( \left| \zeta ( 2\alpha ^{\prime }+\alpha -\frac{1}{2}+i
\left( \frac{1}{2}-\alpha \right) \tan  \theta ) \right| \right) d\theta
\label{XRef-Equation-21205816}
\end{aligned}
\end{equation}

\noindent where $\alpha$, $\alpha ^{\prime }$ are related by (\ref{XRef-Equation-21201724}).

Thus the integral along the vertical line of abscissa $2\alpha ^{\prime
}+\alpha -\frac{1}{2}$ is the same as along the vertical line of
abscissa $\alpha +\frac{1}{2}$ as long as $\frac{1}{1-2\alpha }=\zeta
( 2\alpha ^{\prime }) $ if RH is true and vice versa. We may now
formulate the second Theorem expressing such an equivalence. Notice
that in the above formulas the Lorentz measure is now disappeared
and the ``hidden symmetry'' appears as a global axial symmetry i.e.
which is not pointwise but which is related to an important arithmetical
function given by (\ref{XRef-Equation-21201724}) above.
\begin{theorem}

The RH is equivalent to the existence of a ``hidden symmetry'' given
by:\label{XRef-Theorem-212145}
\[
\varphi ( \alpha ) =\varphi ( \alpha ^{\prime }) 
\]

where
\[
\alpha =\frac{1}{2}\left( 1-\frac{1}{\zeta ( 2\alpha ^{\prime })
}\right) =\frac{1}{2}\left( 1-\prod \limits_{p\,\mathrm{prime}}\left(
1-\frac{1}{p^{2\alpha ^{\prime }}}\right) \right) 
\]

or written in another way:
\[
\alpha =\frac{1}{2}\left( 1-\frac{1}{\zeta ( 2\alpha ^{\prime })
}\right) =\frac{1}{2}\left( 1-\sum \limits_{n=1}^{\infty }\frac{\mu
( n) }{n^{2\alpha ^{\prime }}}\right) 
\]

Moreover
\[
\varphi ( \alpha ) =\ln \left( \frac{1}{1-2\alpha }\right) =\varphi ( \alpha
^{\prime }) =\ln ( \zeta ( 2\alpha ^{\prime }) ) 
\]
\end{theorem}

In (\ref{XRef-Equation-2120580}) and (\ref{XRef-Equation-21205816}),
the Zeta function appears itself, through $\alpha$ and $\alpha ^{\prime
}$ or only $\alpha$, as argument of the Zeta function of complex
argument, by means of the primes. It may be added that a part the
scaling factor $\frac{1}{2}-\alpha $ on the vertical line (the same
outside as well as inside the critical strip) the ``hidden symmetry''
says that for any $0<\alpha <\frac{1}{2}$, the potential of a test
charged vertical filament placed inside the critical strip at the
position $x=\frac{1}{2}+\alpha $ is the same as that of the filament
placed outside the critical strip at the position $x=2\alpha ^{\prime
}+\alpha -\frac{1}{2}$, where $\alpha$, $\alpha ^{\prime }$ are
related by (\ref{XRef-Equation-21201724}). To the best of our knowledge,
the formulation given by Theorem \ref{XRef-Theorem-212145}, even
if equivalent to the RH, is new. It has an electrostatic interpretation
and go more in the direction to believe that the RH my be true.

As an illustration of the kind of convergence in the numerical treatment,
we treat a special case below and add the plots of the two corresponding
periodic functions for a special value of $\alpha$.

\section{Numerical experiment as an illustration of convergence}

We now perform a numerical experiment and consider the case where
$\zeta ( 2\alpha ^{\prime }) =e\cong 2.71828$. Then up to some decimals
we find that $2\alpha ^{\prime }\cong 1.47446$ and so $\alpha =\frac{1}{2}(1-\frac{1}{e})$.
So for this special choice of $\alpha$, the values of the potential
$\varphi$ are the same at $\rho =\alpha +\frac{1}{2}=1-\frac{1}{2e}\cong
0.81606$ (in the critical strip) and at $\rho =2\alpha ^{\prime
}+\alpha +\frac{1}{2}\cong 1.29052$ (outside the critical strip).
We first give the plots of the two functions to integrate in (\ref{XRef-Equation-2120580})
and (\ref{XRef-Equation-21205816}) as a function of $\theta$ (up
to $\theta =2\pi $), to show the periodicity.
\begin{figure}[h]
\begin{center}
\includegraphics{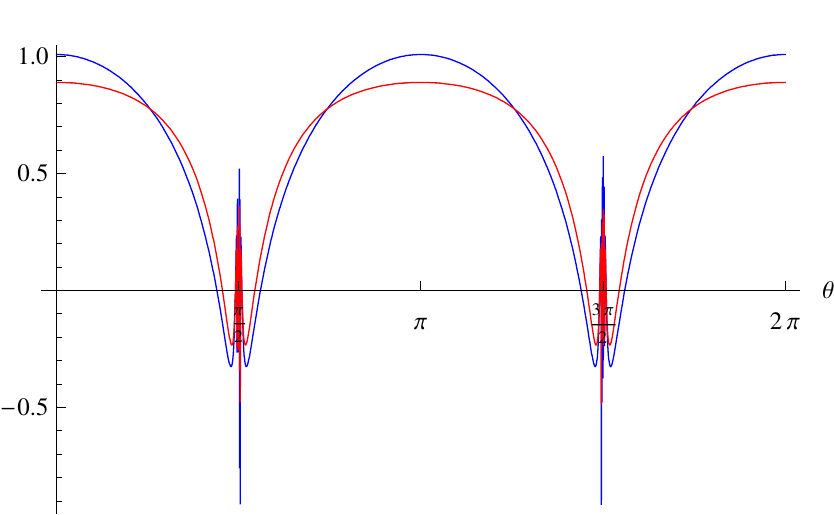}

\end{center}
\caption{The two functions to integrate: the function in (\ref{XRef-Equation-2120580})
[blue] and in (\ref{XRef-Equation-21205816}) [red] }

\end{figure}

Integration of the two functions from $0$ to $0.999\frac{\pi }{2}\cong
1.56923$ gives for the first the value $0.999995$ and for the second
the value $0.999997$. We know the last value should converge unconditionally
to $\ln ( \zeta ( 2\alpha ^{\prime }) ) =\ln ( e) =1$. This illustrate
the kind of convergence involved. Notice that in both cases the
``height'' of integration is given by $(\frac{1}{2}-\alpha )\tan
\theta \cong 117.1$, which corresponds to consider the first 37
non trivial zeros of the Zeta function.
\begin{figure}[h]
\begin{center}
\includegraphics{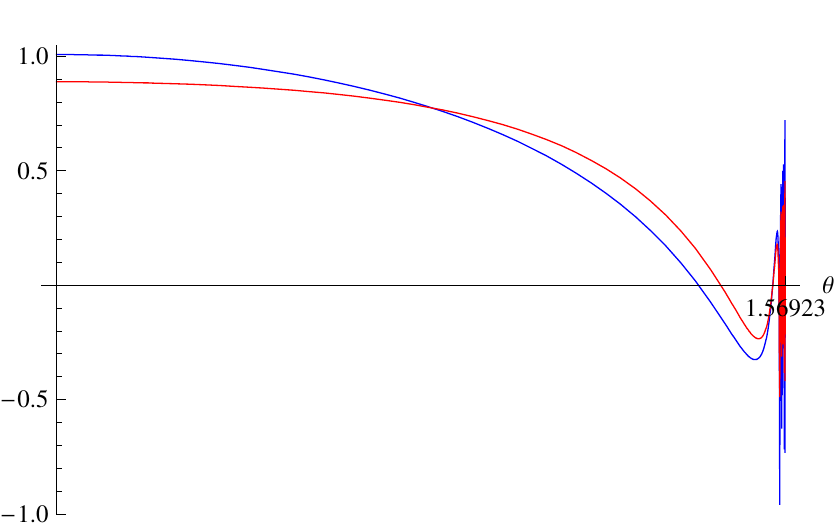}

\end{center}
\caption{Plot of the two functions in the interval $[0,0.999\frac{\pi
}{2}]$}

\end{figure}

Finally, we present the plot of the potentials $\varphi ( \alpha
) =\ln ( \frac{1}{1-2\alpha }) $ and $\varphi ( \alpha ^{\prime
}) =\ln ( \zeta ( 2\alpha ^{\prime }) ) $ which illustrate the divergence
at the right border of the critical strip ($\alpha =\alpha ^{\prime
}=\frac{1}{2}$). We may also consider the ``electrical field'' inside
and outside the critical strip which we define here as:
\begin{equation}
E( \alpha ) =\frac{d}{d \alpha }\varphi ( \alpha ) =\left\{ \begin{array}{lll}
 \frac{2}{1-2\alpha }=2 e ^{\varphi ( \alpha ) } & , & 0<\alpha
<\frac{1}{2} \\
 2\frac{\zeta ^{\prime }( 2\alpha ) }{\zeta ( 2\alpha ) } & , &
\alpha >\frac{1}{2}
\end{array}\right. %
\label{XRef-Equation-2172124}
\end{equation}

Figure 3 summarizes the example we presented.
\begin{figure}[h]
\begin{center}
\includegraphics{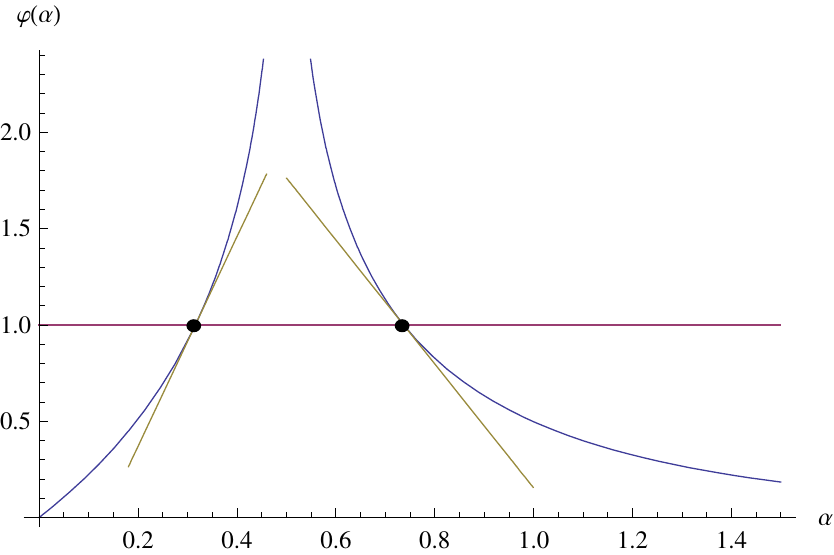}

\end{center}
\caption{The potentials equal to 1 for $\alpha \cong 0.31606$ and
$\alpha \cong 0.737232$ where the slopes of the two tangent straight
lines are given by the magnitude of the corresponding ``electrical
field'' defined by (\ref{XRef-Equation-2172124})}

\end{figure}

Notice that a slightly different definition of the electrical field
i.e. if defined as the gradient of the potential $\varphi$ with
respect to the position $\rho$ (keeping $\rho _{0}$ constant) gives
the following: 
\begin{equation}
E( \alpha ) =\frac{d}{d \rho }\varphi ( \alpha ) =\left\{ \begin{array}{lll}
 \gamma +\frac{1}{1-\rho +\rho _{0}} & , & 0<\alpha <\frac{1}{2}
\\
 \frac{\zeta ^{\prime }( \rho +\rho _{0}) }{\zeta ( \rho +\rho _{0})
} & , & \alpha >\frac{1}{2}
\end{array}\right. 
\end{equation}

Moreover with our choice we obtain:
\[
E( \alpha ) =\left\{ \begin{array}{lll}
 \gamma +\frac{1}{1-2\alpha } & , & 0<\alpha <\frac{1}{2} \\
 \frac{\zeta ^{\prime }( 2\alpha ) }{\zeta ( 2\alpha ) } & , & \alpha
>\frac{1}{2}
\end{array}\right. 
\]

\section{Conclusions}

In this work we have first extended some integral formulas for the
logarithm of the Zeta function\ \ by means of more general Lorentz
measures and obtained two new relations equivalent to the RH. Then
we have introduced a kind of ``hidden symmetry'' which relates the
integrals between two values of the abscissa, one inside the other
outside the critical strip. A simple numerical experiment has been
presented as illustration of the kind of convergence involved.

We note, the truth of such a symmetry is still equivalent to the
truth of the RH, but in the new formulation we have found, there
is the appearance of the primes into the critical strip by means
of the Zeta function calculated outside the critical strip.

The above symmetry is weaker then the stronger Riemann symmetry
\cite{26} which for the Xi\ \ function is given pointwise by $\xi
( s) =\xi ( 1-s) $ for any complex argument {\itshape s} of the
Zeta function. Our weaker symmetry connects the interval [$\frac{1}{2},1$[
in the critical strip with the infinite interval ]$1,\infty $[ outside
the critical strip. The work will be continued with the study of
new integral relations\ \ with a more general class of ``measure''
in the integration of the logarithm of the Zeta function \cite{27}.

\appendix


\begin{thebibliography}{000}
\bibitem{28} F.T. Wang, A note on the Riemann Zeta-Function, Bull.
Amer. Math. Soc., {\bfseries 52}, 1946, 319-321\label{28}
\bibitem{22} V.V. Volchov, On a equality equivalent to the Riemann
Hypothesis, Ukranian Mathematical Journal, {\bfseries 47}, 1995,
422\label{22}
\bibitem{23} M. Balazard, E. Saias, M. Yor, Notes sur la fonction
de Riemann, Advances in Mathematics, {\bfseries 143}, 1995, 422\label{23}
\bibitem{24} D. Merlini, The Riemann Magneton of the Primes, Chaos
and complexity Letters, Nova Science Publishers, New York, Vol.
2, Number 1, 2006, 93\label{24}
\bibitem{25} I.S. Gradshteyn, I.M. Ryzhik, Table of Integrals Series
and Products, Academic Press, New York, 1965, 560\label{25}
\bibitem{26} B. Riemann, Oeuvres Mathematiques, Edition Jacques
Gabay, 1990, 165\label{26}
\bibitem{27} S. Beltraminelli, D. Merlini, S. Sekatskii, in preparation,
2008\label{27}
\end{thebibliography}
\end{document}